\documentclass[twocolumn,amssymb,aps,nofootinbib]{revtex4}
\usepackage{times,mathptmx,amsthm,amsmath}
\usepackage{pstricks,pst-node,graphics}
\begin{document}

\newcommand{\BC}{{1}}
\newcommand{\BP}{{2}}
\newcommand{\EKLP}{{3}}
\newcommand{\Ka}{{4}}
\newcommand{\Ku}{{5}}
\newcommand{\Pro}{{6}}
\newcommand{\RR}{{7}}
\newcommand{\SP}{{8}}
\newcommand{\TF}{{9}}

\newcommand{\R}{{\bf R}}
\newcommand{\Q}{{\bf Q}}
\newcommand{\N}{{\bf N}}
\newcommand{\Z}{{\bf Z}}
\newcommand{\A}{{\cal A}}
\newcommand{\La}{{\lambda}}
\newcommand{\lambdadet}{{\det\,\!_\lambda}}

\begin{center}
Lambda-determinants and domino-tilings \\
$ $ \\
James Propp \\
University of Wisconsin \\
May 16, 2004 
\end{center}

\begin{center}
{\it dedicated to the memory of David Robbins}
\end{center} 

\vspace{0.5in}

{\sc Abstract:} 
Consider the $2n$-by-$2n$ matrix $M=(m_{i,j})_{i,j=1}^{2n}$
with $m_{i,j} = 1$ for $i,j$ satisfying $|2i-2n-1|+|2j-2n-1| \leq 2n$
and $m_{i,j} = 0$ for all other $i,j$,
consisting of a central diamond of 1's 
surrounded by 0's.
When $n \geq 4$, the $\lambda$-determinant of the matrix $M$ 
(as introduced by Robbins and Rumsey [\RR])
is not well-defined.
However, if we replace the $0$'s by $t$'s,
we get a matrix whose $\lambda$-determinant
is well-defined and is a polynomial in $\lambda$ and $t$.
The limit of this polynomial as $t \rightarrow 0$
is a polynomial in $\lambda$ whose value at $\lambda=1$
is the number of domino tilings of a $2n$-by-$2n$ square.

{\large
\bigskip
\noindent
1. Lambda-determinants ...
\bigskip
\normalsize}

In section 5 of their article [\RR],
David Robbins and Howard Rumsey, Jr.\ defined
a generalization of the determinant of a matrix,
which they dubbed the $\lambda$-determinant.
It is a rational function of the entries of the matrix
along with an extra parameter, $\lambda$;
when $\lambda$ is set equal to $-1$,
one obtains the ordinary determinant of the matrix,
at least in the case where all matrix entries are non-zero.

In this article I will consider certain matrices
with many vanishing entries.
For these matrices, one cannot 
apply Robbins and Rumsey's definition literally,
but there is still a natural way
to attempt to compute the $\lambda$-determinant,
by replacing the zeroes by an indeterminate $t$, 
taking the $\lambda$-determinant of the $t$-perturbed matrix,
and then taking the limit as $t \rightarrow 0$.
In particular, I will give a one-parameter family of $0,1$-matrices
whose $n$th member is a $2n$-by-$2n$ matrix whose $\lambda$-determinant,
defined by this continuity method and then specialized to $\lambda = 1$,
is the number of domino-tilings of a $2n$-by-$2n$ square.

Let us start by recalling Robbins and Rumsey's first, recursive definition 
of the $\lambda$-determinant.
If $M$ is a 1-by-1 matrix, its $\lambda$-determinant is
its sole entry.
Now suppose $M$ is an $n$-by-$n$ matrix, with $n>1$.
Let $M_{NW}$, $M_{NE}$, $M_{SW}$, and $M_{SE}$
denote the $\lambda$-determinants of the $(n-1)$-by-$(n-1)$
connected submatrices in the northwest, northeast, southwest, and southeast
corners of $M$,
and let $M_C$ be the $\lambda$-determinant of the central connected
$(n-2)$-by-$(n-2)$ submatrix of $M$
(we take $M_C=1$ in the case $n=2$).
As long as $M_C$ is non-zero, we define the $\lambda$-determinant of $M$ as
\begin{equation}
\lambdadet M 
= (M_{NW} M_{SE} + \lambda M_{NE} M_{SW})/M_C.
\end{equation}
Using this definition, we can calculate
the $\lambda$-determinant of the matrix
$$
\left( \begin{array}{cc}
1 & 1 \\
1 & 1 \end{array} \right)
$$
as $1+\lambda$
and the $\lambda$-determinant of the matrix
$$
\left( \begin{array}{ccc}
1 & 1 & 1 \\
1 & 1 & 1 \\
1 & 1 & 1 \end{array} \right)
$$
as $((1+\lambda)(1+\lambda)+\lambda(1+\lambda)(1+\lambda))/(1) 
= (1+\lambda)^3$.
As was pointed out by Robbins and Rumsey
(and is easy to check by induction),
the $\lambda$-determinant of the $n$-by-$n$ all-ones matrix
is $(1+\lambda)^{n(n+1)/2}$.

When one is using (1) to calculate $\lambda$-determinants,
a subtle distinction becomes important,
namely, the distinction between working in $\Q(\lambda)$
throughout the recursion
and substituting a particular value of $\lambda$ at the end
(on the one hand),
and using that particular value of $\lambda$
when performing the recursion (on the other).
Consider, for instance,
the 4-by-4 matrix whose entries are all 1's.
Its $\lambda$-determinant is $(1+\lambda)^6$,
and if we put $\lambda=-1$, we get 0.
However, if we were to use $\lambda=-1$
in carrying out the recurrence,
we would run into trouble,
since $M_{NW}$, $M_{NE}$, $M_{SW}$, $M_{SE}$, and $M_C$ all vanish;
the $(-1)$-determinant of the matrix
is given by the indeterminate expression $((0)(0)-(0)(0))/0$.
If our goal is to be make sense of the $\lambda$-determinant
over as broad a class of matrices as possible,
clearly we should work in $\Q(\lambda)$ whenever we can.

Robbins and Rumsey give another, non-recursive formula
for the $\lambda$-determinant
of an $n$-by-$n$ matrix $M=(m_{i,j})_{i,j=1}^n$:
\begin{equation}
\lambdadet (M) =
\sum_{B \in \A_n} \lambda^{P(B)} (1+\lambda)^{N(B)} M^B.
\end{equation}
Here $\A_n$ is the set of $n$-by-$n$ alternating-sign matrices
$B=(b_{i,j})_{i,j=1}^n$,
$P(\cdot)$ and $N(\cdot)$ are integer-valued functions on $\A_n$, and
$$M^B = \prod_{i=1}^n \prod_{j=1}^n m_{i,j}^{b_{i,j}}.$$
In more detail:
An alternating-sign matrix is a matrix of $+1$'s, $-1$'s, and $0$'s
such that in each row and column,
the non-zero entries alternate in sign,
beginning and ending with a $+1$ (which may be the same entry).
If $B$ is an $n$-by-$n$ alternating-sign matrix,
we define its inversion number as
$I(B) = \sum b_{i,j} b_{r,s}$
where the sum is over all $1 \leq i,j,r,s \leq n$
with $i<r$ and $j>s$
(if $B$ is a permutation matrix,
this coincides with the ordinary inversion statistic).
We also define $N(B)$ as the number of negative entries in $B$,
and $P(B)$ as $I(B)-N(B)$.

The summation formula (1) for $\lambdadet M$
has exactly the same domain of applicability
as the recursive formula (2), if we work over $\Q(\lambda)$
(and if the matrix entries themselves do not depend on $\lambda$);
specifically, both formulas apply and give the same answer
as long as the entries of $M$
that are in the central $(n-2)$-by-$(n-2)$ submatrix
are all non-zero
(these are precisely the positions in an $n$-by-$n$ alternating-sign matrix
where the entry $-1$ can occur). 
However, if we use specific values of $\lambda$ in the recursion,
the recursive formula can run into problems
where the summation formula does not.
For instance, consider the 4-by-4 matrix with all entries equal to 1;
as we have seen, if we try to compute its $(-1)$-determinant
using the specialization of the recurrence to $\lambda=-1$,
we get an indeterminate result,
whereas the formula in terms of alternating-sign matrices gives 0.

Unfortunately, neither of Robbins and Rumsey's two definitions
works when central entries of $M$ are equal to zero.
Nor is a straightforward appeal to continuity going to help us
to define $\det_{\lambda} M$ for every $M$.  Consider for instance 
the family of matrices
$$
M_c (t) =
\left( \begin{array}{ccc}
t & t & t \\
t & t^4/c & t \\
t & t & t \end{array} \right) .
$$
For $c$ and $t$ non-zero,
$\lambdadet M_c (t) = 
(c\La+c\La^2)+(2\La+2\La^2)t^3+(1/c+\La^3/c)t^6$,
which converges to
$c\La+c\La^2$ as $t \rightarrow 0$.
This limit depends on the value $c$.
Hence, if we attempt to define the $\lambda$-determinant of
the three-by-three all-zeroes matrix
by taking a trajectory through that matrix
in the space of three-by-three matrices
and invoking continuity,
the limit will depend on the trajectory we choose,
and may even fail to exist.

Clearly the principled thing to do would be
to study continuity properties of the
$\lambda$-determinant,
and I hope others will adopt this approach
and undertake a more systematic study
of what happens when different trajectories
through a bad matrix are taken;
this may have some bearing on the issue
of how Dodgson condensation can be applied
to matrices with vanishing connected minors
(i.e., with singular connected submatrices).
However, in this article I will take an easier path
and restrict attention to the trajectory through $M$
in a particular direction.
Specifically, I will replace the zeroes in $M$
by a new variable, $t$,
and see what happens to the $\lambda$-determinant
of the perturbed matrix 
(the ``$t$-perturbation of $M$'')
as $t \rightarrow 0$.
This may be unprincipled, but it is easy to compute.
Moreover, for many matrices $M$,
the resulting rational function in $t$
is not just a Laurent polynomial in $t$
(as it must be, because of the summation formula),
but is in fact an ordinary polynomial in $t$.
In this case, the $\lambda$-determinant of $M$
can be defined as the constant term of this polynomial.

For example, consider the eight-by-eight matrix
$$
\left( \begin{array}{cccccccc}
0&0&0&1&1&0&0&0 \\
0&0&1&1&1&1&0&0 \\
0&1&1&1&1&1&1&0 \\
1&1&1&1&1&1&1&1 \\
1&1&1&1&1&1&1&1 \\
0&1&1&1&1&1&1&0 \\
0&0&1&1&1&1&0&0 \\
0&0&0&1&1&0&0&0
\end{array} \right) .
$$
It has interior zeroes,
so its $\lambda$-determinant cannot be computed
by the original formulas (1) and (2).
However, if we replace all 0's by $t$,
we get a polynomial in $\lambda$ and $t$ with 191 terms.
Replacing $t$ by 0,
we get a polynomial in $\lambda$ with a mere 17 terms,
which evaluates to 12988816
when we set $\lambda=1$.

{\large
\bigskip
\noindent
2. ... and Domino-Tilings
\bigskip
\normalsize}

12988816 is also the number of ways to cover an 8-by-8 square
with 32 1-by-2 rectangles, commonly known as dominos.
More generally, the number of ways to cover a $2n$-by-$2n$ square
with $2n^2$ dominos was computed by Temperley and Fisher [\TF]
and (simultaneously) by Kasteleyn [\Ka],
and is given by the double product
$$\prod_{j=1}^{n} \prod_{k=1}^{n}
\left(4 \cos^2 \frac{\pi j}{2n+1} + 4 \cos^2 \frac{\pi k}{2n+1}\right).$$

It turns out that if one considers the sub-region 
of the $2n$-by-$2n$ square that consists of 
the 2 central cell in the first and last rows,
the 4 central cells in the second and second-from-last rows,
the 6 central cells in the third and third-from-last rows,
etc.,
one gets a region whose domino-tilings are enumerated
by a much simpler expression, namely 
$$2^{n(n+1)/2}.$$
This region is called the Aztec diamond of order 8,
and was first studied in detail in [\EKLP],
although earlier occurrences of the shape appear in the literature.
As an example, consider the case $n=4$.
If we associate the 64 cells of the 8-by-8 square
with the entries of the 8-by-8 matrix
considered at the end of the previous section,
then the entries that contain 1's correspond to the cells
that belong to the Aztec diamond of order 4.

It turns out that the domino-tilings of the Aztec diamond of order $n$
are intimately related to the $\lambda$-determinant 
of the generic $n$-by-$n$ matrix.
(If the meaning of ``generic'' is unclear,
then the reader should imagine that we are working
over the ring $\Q(\lambda,a,b,c,\dots)$,
where $a,b,c,\dots$ are the entries of the matrix;
then the entries of the matrix are perforce non-zero,
so the summation formula for the $\lambda$-determinant 
is unproblematical in this context.)
Specifically, if we apply the distributive rule
to the Robbins-Rumsey summation formula for $\det_{\lambda}$,
and ignore the fact that multiplication is commutative,
each term
$\lambda^{P(B)} (1+\lambda)^{N(B)} M^B$
expands to a sum of $2^{N(B)}$ monomials
of the form $\lambda^k M^B$.
The resulting sum has $2^{n(n+1)/2}$ terms,
and these terms can be put into 1-1 correspondence
with the $2^{n(n+1)/2}$ domino-tilings
of the Aztec diamond of order $n$.
(See [\EKLP] for details.)

If it seems paradoxical that the $\lambda$-determinant
of a square matrix of 1's counts domino-tilings of an Aztec diamond,
whereas the $\lambda$-determinant of the 0,1-matrix
whose 1's form an Aztec diamond counts domino-tilings of a square,
it may be helpful to imagine rotating the Aztec diamond by 45 degrees.
Here is one of the 64 domino-tilings of the Aztec diamond of order 3:

$$
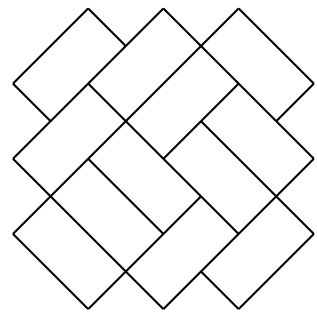
\pspicture(-2.2,-2.2)(2.2,2.2)
\psset{xunit=.5cm,yunit=.5cm}
\psline(-4,2)(-2,4)
\psline(-4,0)(0,4)
\psline(-4,-2)(2,4)
\psline(0,0)(2,2)
\psline(-2,-4)(1,-1)
\psline(3,1)(4,2)
\psline(0,-4)(4,0)
\psline(2,-4)(4,-2)
\psline(-4,-2)(-2,-4)
\psline(-4,0)(0,-4)
\psline(-4,2)(-3,1)
\psline(-2,0)(0,-2)
\psline(1,-3)(2,-4)
\psline(-2,2)(2,-2)
\psline(-2,4)(-1,3)
\psline(1,1)(4,-2)
\psline(0,4)(4,0)
\psline(2,4)(4,2)
\endpspicture
$$

\noindent
We can replace the tiling problem
by the dual matching problem.
We define an Aztec diamond graph whose vertices
correspond to the cells of the Aztec diamond,
with an edge between two vertices when
the two corresponding cells are adjacent.
Then a tiling of the Aztec diamond of order $n$
corresponds to a perfect matching of the Aztec diamond graph,
that is, a collection of edges with the property that
each vertex of the graph belongs to exactly one
of the edges in the collection.

To illustrate, here is the Aztec diamond graph of order 3:

$$
\pspicture(-1.8,-1.8)(1.8,1.8)
\psset{xunit=.5cm,yunit=.5cm}
\pspolygon(-3,-2)(2,3)(3,2)(-2,-3)
\pspolygon(-3,0)(0,3)(3,0)(0,-3)
\pspolygon(-3,2)(-2,3)(3,-2)(2,-3)
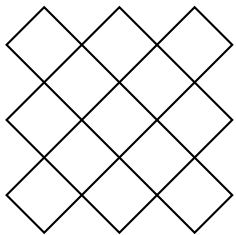
\endpspicture
$$

\noindent
And here is the perfect matching of the Aztec diamond graph of order 3
that corresponds to the domino-tiling of the Aztec diamond shown in
the first figure:

$$
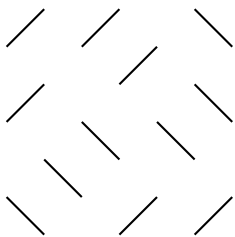
\pspicture(-1.8,-1.8)(1.8,1.8)
\psset{xunit=.5cm,yunit=.5cm}
\psline(-3,2)(-2,3)
\psline(-3,-2)(-2,-3) 
\psline(3,2)(2,3)
\psline(3,-2)(2,-3)
\psline(-3,0)(-2,1)
\psline(0,-3)(1,-2)
\psline(-1,-2)(-2,-1)
\psline(0,-1)(-1,0)
\psline(-1,2)(0,3)
\psline(0,1)(1,2)
\psline(2,-1)(1,0)
\psline(3,0)(2,1)
\endpspicture
$$

\noindent
Looking back at the original figure, note that if the
northwest, northeast, southwest, and southeast tiles are removed,
what's left is a domino-tiling of the 4-by-4 square.
More generally, there is a one-to-one correspondence
between domino-tilings of the $2n$-by-$2n$ square
and domino-tilings of the Aztec diamond of order $2n$
in which there are $n(n-1)/2$ forced tiles
in each of the four corners
(slanting from southwest to northeast
in the northwest and southeast corners,
and slanting from northwest to southeast 
in the southwest and northeast corners).
The latter in turn correspond to perfect matchings
of the Aztec diamond graph of order $2n$
in which there are $n(n-1)/2$ forced edges in each corner.

One way we can force some edges to be present
is to force other edges to be absent,
and one way to force edges to be absent
is to work in the setting of weighted enumeration.
Given an assignment of non-negative weights
to the edges of a graph,
we define the weight of an individual perfect matching
as the product of the weights of its edges.
In the article [\Pro],
I showed how the method of [\EKLP]
could be adapted to the general problem
of finding the sum of the weights
of all the perfect matchings of an
edge-weighted Aztec diamond graph.
I applied this to the case of $2n$-by-$2n$ squares,
setting some edge-weights equal to 1
and the rest equal to 0,
in the fashion shown below
for an Aztec diamond graph of order 5
(the edges shown get weight 1,
and the rest get weight 0):

$$
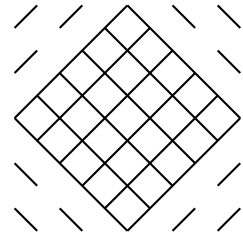
\pspicture(-1.8,-1.8)(1.8,1.8)
\psset{xunit=.3cm,yunit=.3cm}
\psline(-5,0)(0,5)
\psline(-4,-1)(1,4)
\psline(-3,-2)(2,3)
\psline(-2,-3)(3,2)
\psline(-1,-4)(4,1)
\psline(0,-5)(5,0)
\psline(-5,0)(0,-5)
\psline(-4,1)(1,-4)
\psline(-3,2)(2,-3)
\psline(-2,3)(3,-2)
\psline(-1,4)(4,-1)
\psline(0,5)(5,0)
\psline(-5,-4)(-4,-5)
\psline(-5,-2)(-4,-3)
\psline(-3,-4)(-2,-5)
\psline(5,-4)(4,-5)
\psline(5,-2)(4,-3)
\psline(3,-4)(2,-5)
\psline(-5,4)(-4,5)
\psline(-5,2)(-4,3)
\psline(-3,4)(-2,5)
\psline(5,4)(4,5)
\psline(5,2)(4,3)
\psline(3,4)(2,5)
\endpspicture
$$

\noindent
However, what was missing from that account
was an explanation of the link with
$\lambda$-determinants.

A good way to understand this link,
without using all the machinery of generalized domino-shuffling,
is to directly rely on Kuo's method of graphical condensation [\Ku].
Let $G$ be a weighted Aztec diamond graph of order $n$.
Let $G_{NW}$ be the weighted Aztec diamond graph of order $n-1$
derived from $G$ by elminating the southernmost $n$ vertices,
the southernmost $2n$ edges, the easternmost $n$ vertices, and
the easternmost $2n$ edges.  Define $G_{NE}$, $G_{SW}$, and $G_{SE}$
analogously.  Define $G_C$ to be the weighted Aztec diamond graph
of order $n-2$ derived from $G$ by eliminating all of the aforementioned
vertices and edges.  Lastly, define $w_{NW}$ to be the weight of the
northwesternmost edge, and define $w_{NE}$, $w_{SW}$, and $w_{SE}$
analogously.  Then Kuo's formula can be stated as
\[
\begin{array}{l}
W(G) W(G_C) = \\
w_{NE} w_{SW} W(G_{NW}) W(G_{SE}) + w_{NW} w_{SE} W(G_{NE}) W(G_{SW}),
\end{array}
\]
where $W(\cdot)$ denotes the sum of the weights of the perfect matchings
of the graph in question.
In our application, the edge-weights
$w_{NW}$, $w_{NE}$, $w_{SW}$, and $w_{SE}$
are all 0's and 1's.
This formula becomes a recurrence relation
if one divides both sides by $W(G_C)$
(though this is only sensible if $W(G_C)$ is non-zero).
If one runs the recurrence for the case at hand,
the prefactors $w_{NE} w_{SW}$ and $w_{NW} w_{SE}$
are always equal to 1.
Hence the recurrence takes the simplified form
$$ W(G) = (W(G_{NW} W(G_{SE}) + W(G_{NE}) W(G_{SW})) / W(G_C), $$
which is equation (1) in the special case $\lambda=1$.
This gives us a combinatorial proof
of the main claim of this paper,
namely, that the number of domino-tilings of the $2n$-by-$2n$ square
can be calculated as the $(+1)$-determinant
of the $2n$-by-$2n$ matrix
whose $i,j$th entry (for $1 \leq i,j \leq 2n$) 
is 1 if $|2i-2n-1|+|2j-2n-1| \leq 2n$ and is 0 otherwise.

Kuo's method also gives us a combinatorial interpretation 
to all the numbers that occur in the course of evaluating
the $(+1)$-determinant of the matrix (by applying (1) with
$\lambda=1$): these numbers count domino-tilings of regions
obtained from the Aztec diamond by eliminating the cells 
that lie in certain bands bordering the boundaries of the 
region.  

For instance, consider the process of computing the 
$(+1)$-determinant of the 4-by-4 matrix
$$
M= 
\left( \begin{array}{cccc}
0&1&1&0 \\
1&1&1&1 \\
1&1&1&1 \\
0&1&1&0
\end{array} \right) .
$$
The quantities that turn up
in the recursive application of (1)
are precisely the $(+1)$-determinants 
of all the connected submatrices of $M$.
$M$ itself does double-duty as the matrix of
$(+1)$-determinants of all the 1-by-1 submatrices of $M$.
The $(+1)$-determinants of the connected 2-by-2 submatrices of $M$
form the 3-by-3 matrix
$$
\left( \begin{array}{ccc}
1&2&1 \\
2&2&2 \\
1&2&1 
\end{array} \right)
$$
whose entries are given by the $\lambda$-condensation formula (1)
in the special case $\lambda=+1$.
Turning the crank again gives the 2-by-2 matrix
$$
\left( \begin{array}{cc}
6&6 \\
6&6
\end{array} \right)
$$
and turning it a final time gives the 1-by-1 matrix
$$
\left( \begin{array}{c}
36 
\end{array} \right)
$$
whose sole entry is the number of tilings of
the 4-by-4 square.
The 6's count the domino-tilings of the region
obtained from the 4-by-4 square
by removing three cells from one corner
and three cells from an adjoining corner.

This combinatorial interpretation of the
$(+1)$-determinants of the connected submatrices of $M$
makes it clear that Robbins and Rumsey's original recursive formula
for the $\lambda$-determinant can be applied to all of the matrices
in our one-parameter family, provided we interpret the indeterminate 
expression $0/0$ as 0 whenever it crops up; each such occurrence
corresponds to a subgraph of the weighted Aztec diamond graph
whose perfect matchings all have weight 0.

The one thing missing from this explanation
is an explicit discussion of
the behavior of the $(+1)$-determinant
of the $t$-perturbed matrix.
We need to know that 
it is a polynomial in $t$.
But this follows from the fact that it is equal to
the sum of the weights of all the perfect matchings
of the weighted graph.
Indeed, all the rational functions of $t$
that occur during the recursion
are polynomials in $t$,
for the same reason.

This takes care of the case $\lambda=+1$:
the inclusion of $t$ has solved
the indeterminacy problem.
It follows a fortiori that for generic $\lambda$,
inclusion of $t$ also lets us carry out the recurrence (1) 
without encountering indeterminacy.
(However, if one puts $\lambda=-1$,
one still encounters indeterminacy,
on account of the cancellations that occur.)

One consequence of the main result of this paper
is that if one takes any $2n$-by-$2n$ alternating-sign matrix 
$A=(a_{i,j})_{i,j=1}^{2n}$
and sums those entries $a_{i,j}$ for which $|2i-2n-1|+|2j-2n-1| \leq 2n$,
one gets a non-negative sum.
A similar claim holds for odd-by-odd alternating-sign matrices:
if one takes any alternating-sign matrix 
$A=(a_{i,j})_{i,j=1}^{2n+1}$
and sums those entries $a_{i,j}$ for which $|i-n-1|+|j-n-1| \leq n$,
one gets a non-negative sum.
This allows one to use the $t$-perturbation trick
to extend the definition of the $\lambda$-determinant
to the odd-by-odd matrix $M=(m_{i,j})_{i,j=1}^{2n+1}$
with $m_{i,j}=1$ for those $i,j$ for which $|i-n-1|+|j-n-1| \leq n$
and $m_{i,j}=0$ for all other $i,j$.
(In fact, the $(+1)$-determinant of this matrix
is also equal to the number of domino-tilings
of the $2n$-by-$2n$ square.
This was proved by Trevor Bass and Kezia Charles [\BC].)

Moreover, if one takes the intersection of 
either the odd-by-odd or even-by-even diamond pattern
with any connected square submatrix of $M$,
one gets a smaller pattern which also has the property
that sum of the corresponding entries of an alternating-sign matrix
must be non-negative.
It would be interesting to have a classification of those
partial sums of the entries of an $n$-by-$n$ matrix
that are non-negative for all choices of an alternating-sign matrix.
It would also be desirable to extend the results of this paper
to the context of the octahedron recurrence
with general initial conditions, as considered in [\SP].

\vspace{0.4in}

{\large
\bigskip
\noindent
References
\medskip
\normalsize}

\noindent
{\bf [\BC]} \ \ T.\ Bass and K.\ Charles, 
unpublished memo;
{\tt www.math.wisc.edu/$\sim$propp/reach/Charles/}
{\tt octagon.pdf}.

\medskip
\noindent
{\bf [\BP]} \ \ D.\ Bressoud and J.\ Propp,
How the alternating sign matrix conjecture was solved,
\emph{Notices of the AMS} {\bf 46} (1999), 637--646;
{\tt www.ams.org/notices/199906/fea-bressoud.pdf}.

\medskip
\noindent
{\bf [\EKLP]} \ \ 
N.\ Elkies, G.\ Kuperberg, M.\ Larsen, and J.\ Propp, 
Alternating sign matrices and domino tilings, 
{\it J.\ Algebraic Combin.\/} {\bf 1} (1992), 111--132 and 219--234;
{\tt www.math.wisc.edu/$\sim$propp/aztec.ps.gz}.

\medskip
\noindent
{\bf [\Ka]} \ \ P.~W.\ Kasteleyn, 
The statistics of dimers on a lattice, I. 
The number of dimer arrangements on a quadratic lattice,
{\it Physica} {\bf 27} (1961), 1209--1225.

\medskip
\noindent
{\bf [\Ku]} \ \ E.\ Kuo,
Applications of graphical condensation for enumerating matchings 
and tilings, arXiv: {\tt math.CO/0304090}; to appear in
{\it Theoret.\ Comp.\ Sci.\/} {\bf 319} (2004), 29--57.

\medskip
\noindent
{\bf [\Pro]} \ \ J.\ Propp, Generalized domino-shuffling,
{\it Theoret.\ Comp.\ Sci.\/} {\bf 303} (2003), 267--301.

\medskip
\noindent
{\bf [\RR]} \ \ D.~P.\ Robbins and H.\ Rumsey, Jr.,
Determinants and alternating sign matrices,
{\it Adv.\ in Math.\/} {\bf 62} (1986), 169--184.

\medskip
\noindent
{\bf [\SP]} \ \ D.\ Speyer,
Perfect matchings and the octahedron recurrence,
arXiv: {\tt math.CO/0403417}.

\medskip
\noindent
{\bf [\TF]} \ \ H.~N.~V.\ Temperley and M.~E.\ Fisher,
Dimer problem in statistical mechanics --- an exact result,
{\it Phil.\ Mag.\/} {\bf 6} (1961), 1061--1063.

\end{document}